\documentstyle[epsf]{article}

\newtheorem{definition}{Definition}
\newtheorem{theorem}{Theorem}
\newtheorem{remark}{Remark}
\newtheorem{lemma}{Lemma}
\newtheorem{proposition}{Proposition}
\newtheorem{corollary}{Corollary}
\newtheorem{example}{Example}

\def\La{\Lambda}
\def\la{\lambda}
\def\Z{\:\mbox{\sf Z} \hspace{-0.82em} \mbox{\sf Z}\,}
\def\Zs{\mbox{\scriptsize \sf Z}  \! \! \mbox{\scriptsize \sf Z}}
\def\Znn{\Z_{\ge0}}
\def\C{\rule[0.5pt]{0.2mm}{6pt}{\hspace{-3.6pt}{\rm C}}}
\def\ie{{\em i.e., }}

\def\P{{\cal P}_L}
\def\G{{\cal G}_L}
\def\Pset{\P(\La_i+\La_j,\La_k)}
\def\Gset{\G(2 \La_0,\La_k)}
\def\Mult#1#2#3{\left[#1 \atop #2 \right]_{#3}}
\def\Mults#1#2#3{\left[{\textstyle {#1 \atop #2} } \right]_{#3}}
\def\MMult#1#2#3{\left[\left[#1 \atop #2 \right]\right]_{#3}}

\def\e{\mbox{e}}

\def\case#1#2{{\textstyle{#1\over #2}}}

\begin{document}

\title{A proof of polynomial identities of type\\
       $\widehat{sl(n)_1} \otimes \widehat{sl(n)_1} /
       \widehat{sl(n)_2}$ }

\author{Omar Foda$^1$, Masato Okado$^2$ and S. Ole Warnaar$^1$}

\date{ \it $^1$Department of Mathematics,\\
       \it The University of Melbourne\\
       \it Parkville, Victoria 3052, Australia\\
       \vskip 4mm
       \it $^2$Department of Mathematical Sciences,\\
       \it Faculty of Engineering Science,\\
       \it Osaka University, Toyonaka,\\
       \it Osaka 560, Japan \\
       \vskip 4mm
       \it Dedicated to the memory of Claude Itzykson.}

\maketitle

\begin{abstract}
We present a proof of polynomial identities related to finite
analogues of the branching functions of the coset
$\widehat{sl(n)_1} \otimes \widehat{sl(n)_1} / \widehat{sl(n)_2}$.
\end{abstract}

\section{Introduction}

Consider the affine algebra $\widehat{sl(n)_{\ell}}$, where
$(n-1)$ is the rank and $\ell$ is the level\footnote{The
standard reference on affine algebras is \cite{Kac}. For
an elementary introduction, see \cite{KMPS}. For reviews and
original references, see \cite{GO}, and references therein.}.
Following \cite{DJMO,JMO}, the branching functions of the coset
\begin{equation}
{\cal C}_{n, \ell_1, \ell_2} =
\widehat{sl(n)}_{\ell_1} \otimes
\widehat{sl(n)}_{\ell_2} /
\widehat{sl(n)}_{\ell_1 + \ell_2}
\end{equation}
are characters of the highest weight modules (HWM's) of $W_n$
algebras \cite{FZ}, where $W_2$ is the Virasoro algebra~\cite{BPZ}.
We are interested in computing these branching functions.

\subsection{$q$-series identities}
An important observation, made independently in \cite{LP} in
the context of affine algebras, and in \cite{StonyBrook} in
the context of branching functions, is that different approaches
to computing the characters lead to completely different
expressions for them. Equating different expressions of the
same character leads to generalizations of the Rogers-Ramanujan
identities. In the present work, we are interested in the
identities related to the branching functions.

\subsubsection{Boson-fermion identities}
Because one side of these identities is generated
using operators that obey bosonic commutation relations,
while the other is generated using operators
that obey fermion-like exclusion principles, these
identities are also known as boson-fermion identities.

In \cite{JMO}, the branching functions of the coset
${\cal C}_{n, \ell_1, \ell_2}$ were obtained by counting certain
configurations, known as {\em weighted paths}. These paths
appear naturally in using the corner transfer matrix
method to solve statistical mechanical models
\cite{BaxterBook}. The expressions obtained are of the
bosonic type. In the present work, we restrict our
attention to the coset ${\cal C}_{n, 1, 1}$, and obtain
expressions for the branching functions by counting the
Ferrers graphs that appear in the crystal base
description of the HWM's of $\widehat{sl(n)_1}$. The
expressions obtained are of the fermionic type, and
finitize the Lepowsky and Primc character formulae~\cite{LP}.

\subsubsection{Polynomial identities}
In fact, we do not work directly in terms of the
characters, which are formal infinite series. Instead, we
work in terms of polynomials which depend on a parameter $L$,
and reduce to the characters in the limit $L \rightarrow
\infty$. In that sense, the identities we obtain are
stronger than identities between characters.

Equating the expressions of \cite{JMO}, and those obtained
in the present work, we obtain polynomial identities; one
for each branching function of ${\cal C}_{n, 1, 1}$. For fixed
$n$, there are ${\cal O}(n^2)$ such functions, and corresponding
identities.
These polynomial identities are generalizations of those
considered by Schur in his approach to proving $q$-series
identities \cite{Schur}.

\subsection{Two ways to count}
Though the Ferrers graphs that we count are in one-to-one
correspondence with the weighted paths, the expressions
that we obtain are different from those of \cite{JMO} because
our approach to counting these objects is inherently different.
We wish to outline the usual method of counting, in order to
emphasize the contrast to ours.

\subsubsection{Indirect counting: Sieving}
In \cite{JMO}, the counting was achieved using a {\em sieving
method} to obtain recurrence relations which can solved. The main
idea of the sieving approach can be summarized as
follows:\footnote{For a detailed introduction to the
sieving method, see \cite{AndrewsBook}.}

Suppose one wishes to count the number of objects in a certain
class $P_0$ which satisfy certain conditions\footnote{For
simplicity, we will refer to a class of objects,
and to its cardinality, by the same symbol. Furthermore, by
\lq\lq evaluate a class of objects\rq\rq, we actually mean
compute its cardinality.}.
This is typically a difficult problem, since the conditions
satisfied by $P_0$ can be quite complicated. However, one
can approach it {\em indirectly} as follows:

As a first step, one considers a larger class of objects $Q_0$,
that includes $P_0$, but satisfies weaker conditions, and hence
is easier to evaluate. Suppose one manages to do that, the next
step would be to evaluate the difference $P_1 = Q_0 - P_0$, and
subtract it to obtain $P_0 = Q_0 - P_1$ (hence the name {\em sieving}).
But evaluating $P_1$ {\em directly} is once again typically
just as hard as the initial problem of evaluating $P_0$.
Hence, it should also be evaluated in two steps: We consider
a larger class of objects $Q_1$ that is
easier to evaluate, and subtract that of the difference
$P_2 = Q_1 - P_1$. We obtain $P_0 = Q_0 - Q_1 + P_2$. It is
easy to see how the above procedure generalizes to give
$P_0 = Q_0 - Q_1 + \cdots + Q_{even} - Q_{odd} + \cdots$

The objects we are interested in--Ferrers graphs and paths
--have dimensions. For larger $i$, $P_i$ typically
contains larger objects. If there are no restrictions on the dimensions
of the objects being counted, then the above sieving procedure continues
indefinitely. If there are such restrictions, then for
sufficiently large $i$, the procedure terminates. Either way,
the procedure amounts to writing a recurrence relation for the set
$\{P_0,P_1,\ldots \}$ and solving it.

\subsubsection{Direct counting: sectoring}
In contrast to the above, the approach used in this paper relies
on a direct counting of the objects of interest. The main idea
is to divide the set of all objects into sectors, each of which is
easier to compute, and then to sum over all sectors. An outline of
this approach is given below.

\subsection{Outline of proof}
\begin{enumerate}
\item Given the set of graphs we wish to count, we
propose to distinguish a certain subset to be called
{\em parent graphs}. The remaining graphs are
called {\em non-parents}.

\item We propose a set of rules which reduces {\em any}
non-parent graph uniquely to a parent graph by removing
nodes from it. Using these rules we can decompose any
non-parent graph into a parent graph plus a set of objects
called {\em $g$-components}.
The rules are such that a parent graph cannot be further
reduced to another parent graph.

\item We show that the above set of rules is invertible.
Each non-parent can be uniquely obtained from a parent
by attaching g-components.
Consequently, the set of non-parents which reduce to a given
parent may be regarded as the {\em descendants} of that parent.

\item From the above, we classify the set of all graphs
into sectors. Each sector contains precisely one parent
plus its descendants.

\item We show that, given a parent graph, the set of all
its descendants is generated by a product over Gaussian
polynomials.

\item Since we know the explicit expression for the
Gaussian polynomials in each sector, summing over all
sectors, with the proper weighting which follows from the
weight of the parent graph, we obtain the desired
generating function of the graphs.
\end{enumerate}

\subsection{Plan of paper}
In \S 2, we outline a number of technical details related
to weighted paths on the set of dominant integral
weights of $\widehat{sl(n)}_2$, and recall the
bosonic generating function as evaluated in \cite{JMO}.
In \S 3, we introduce the main objects of this paper:
K-graphs, and discuss their properties.
In \S 4, we describe the special set of K-graphs called parents.
In \S 5, we describe the graph components to be
added to a parent to generate more general K-graphs,
called descendants.
In \S 6, we describe how
the descendants are obtained from their parent, and
why each graph is either a parent, or descends from
a uniquely-defined one.
In \S 7, we evaluate the number of descendants of a certain
parent.
In \S 8, we obtain
fermionic expressions for the finite analogues of all
branching functions of the coset
${\cal C}_{n, 1, 1}$.
In \S 9 we summarize our results to obtain the main theorem
of this paper: polynomial identities for the finite analogues
of the branching functions. This section also contains a
discussion of our results.
\section{Paths}
In this section, we consider weighted paths on the set of
level-2 dominant integral weights of $\widehat{sl(n)}$, and
recall their generating function as computed in \cite{JMO}.

\subsection{Roots and weights}
We start with some definitions from the theory of affine
algebras \cite{Kac}. Let $\La_i,\alpha_i$ ($i=0,\cdots,n-1$)
and $\delta$, be the fundamental weights, the simple roots,
and the null root of the affine Lie algebra $\widehat{sl(n)}$,
respectively. The subscript $i$ of $\La_i$ can be extended to
$i\in\Z$ by setting $\La_i=\La_{i'}$ for $i\equiv i' (\bmod~n)$.
Let $\widehat{i} = \La_{i+1} - \La_i (i = 0, \cdots, n-1)$
be the weights of the vector representation of $sl(n)$, and
$\rho = \sum_{i=0}^{n-1} \La_i$ be the Weyl vector.

\begin{remark}
For the rest of this work, we will simply use
$a\equiv b$ to indicate $a \equiv b\;(\bmod\;n)$.
\end{remark}

Let $P=\Z\La_0\oplus\cdots\oplus\Z\La_{n-1}\oplus\Z\delta$ be
the weight lattice \cite{Kac,KMPS}. There is an invariant
bilinear form $(\cdot\mid\cdot)$ on $P$ defined by
\begin{equation}
(\La_i\mid\La_j)=\min(i,j)-\frac{i~j}{n},\quad
(\La_i\mid\delta)=1,\quad
(\delta\mid\delta)=0,
\end{equation}
for $0\le i,j\le n-1$.

\begin{figure}[hbt]
\centerline{\epsffile{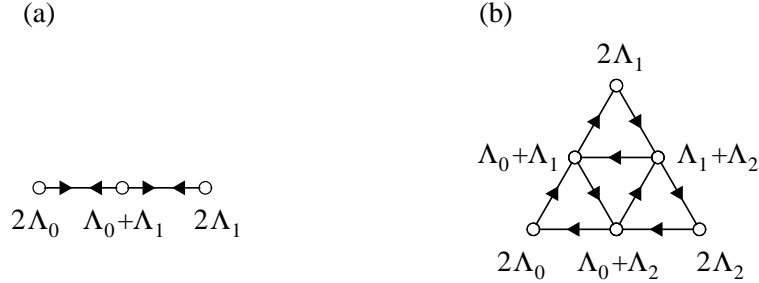}}
\caption{Examples of the set $P_2^+$.
A directed bond from
$\la$ to $\la'$ ($\la,\la'\in P_2^+$) indicates
that a path can go from $\la$ to $\la'$.
(a) $n=2$  (b) $n=3$.}
\label{fig1sec2}
\end{figure}

We are not interested in the full weight lattice, but in
certain restrictions of it:

\begin{definition}[$P_2^+$]
$P_2^+$ is the set of level-2 dominant integral weights, \ie
$P_2^+=\{\La_i+\La_j\mid 0\le i\le j\le n-1\}$.
\end{definition}
Examples of $P_2^+$ in the case of $n=2$, and $3$ are shown in
Figure~\ref{fig1sec2}.

We can define paths on $P_2^+$ as follows:
\begin{definition}[paths]
For $L\in\Znn$, we define a path $p$ as $p=(\la_0,\cdots,\la_L)$
with all $\la_i\in P_2^+$ and $\la_{i+1}-\la_i \in
\{ \widehat{0}, \widehat{1}, \cdots, \widehat{n-1} \}$.
\end{definition}

We are interested in particular sets of paths  of length $L$
defined by
\begin{definition}[$\Pset$]
\begin{equation}
\Pset=\{p=(\la_0,\cdots,\la_L)\mid
\la_0=\La_i+\La_j,\la_L=\La_k+\La_{i+j-k+L}\}.
\end{equation}
For a path $p\in \Pset$ we call
$\La_i+\La_j,\La_k$ and $L$ its {\it initial point},
{\it boundary} and {\it length}, respectively.
\end{definition}
We note that $\Pset$ is a finite analogue (length $L$)
of the set of $(\La_k,\La_{i+j-k})$-restricted paths
of \cite{JMMO,DJO}.

With the  paths in $\Pset$ we associate a special
path $\bar p$ called  the ground-state path,
as follows:
\begin{definition}[ground-state path $\bar{p}$ ]\label{gsp}
$$
\bar{p}=(\La_k+\La_{i+j-k},\La_k+\La_{i+j-k+1},\cdots,
\La_k+\La_{i+j-k+L}) \in\P(\La_k+\La_{i+j-k},\La_k).
$$
\end{definition}
Note that the initial point of the ground-state path
may be different from  that of the paths in $\Pset$.

We can encode a path in terms of a sequence of integers as
follows:
\begin{definition}[sequence of integers]\label{defsoi}
For a path $p=(\la_0,\cdots,\la_L)\in \Pset$ we define
a sequence of integers $\iota(p)=(\mu_0,\cdots,\mu_L)$,
where $\widehat{\mu}_\ell = \la_{\ell+1}-\la_\ell$, and where we
have used $\la_{L+1}=\La_k+\La_{L+i+j-k+1}$. We denote
the element $\mu_\ell$ of $\iota(p)$ by $\iota(p)_\ell$.
\end{definition}
Note that $\iota(\bar{p})$ of $\bar{p}$ in Definition~\ref{gsp}
is given by $\iota(\bar{p})_\ell
\equiv i+j-k+\ell$.
\par\noindent
\begin{example}
The ground state path $\bar{p}$ associated to
${\cal P}_6(\La_i+\La_{-i},\La_0)$ for $n=3$.
\begin{eqnarray*}
\bar{p}&=&(2\La_0,\La_0+\La_1,\La_0+\La_2, 2\La_0,\La_0+\La_1,
\La_0+\La_2,2\La_0) \\
\iota(\bar{p})&=&(0,1,2,0,1,2,0)
\end{eqnarray*}
\end{example}
\par\noindent
\begin{example}
A path in $p^{(1)}\in{\cal P}_6(2\La_0,\La_1)$ for $n=3$.
\begin{eqnarray*}
p^{(1)}&=&(2\La_0,\La_0+\La_1,\La_0+\La_2,\La_1+\La_2,\La_0+\La_1,
          \La_0+\La_2,\La_1+\La_2) \\
\iota(p^{(1)})&=&(0,1,0,2,1,0,2)
\end{eqnarray*}
\end{example}
\begin{example}\label{ex1}
A path in $p^{(2)}\in{\cal P}_6(2\La_0,\La_0)$ for $n=4$.
\begin{eqnarray*}
p^{(2)}&=&(2\La_0,\La_0+\La_1,2\La_1,\La_1+\La_2,2\La_2,
          \La_2+\La_3,\La_0+\La_2) \\
\iota(p^{(2)})&=&(0,0,1,1,2,3,2)
\end{eqnarray*}
\end{example}

\subsection{Weighted paths}
Let $p$ be a path and $\bar{p}$ the ground-state
path associated to $p$,
with integer sequences $\iota(p)=(\mu_0,\ldots,\mu_L)$
and $\iota(\bar{p})=(\bar{\mu}_0,\ldots,\bar{\mu}_L)$, respectively.
We define an energy function $E$ by
\begin{definition}[energy of a path]
\begin{eqnarray}
E(p)&=&
\sum_{\ell=1}^L
\ell\,\big(\theta(\mu_{\ell-1}-\mu_{\ell})-
     \theta(\bar{\mu}_{\ell-1}-\bar{\mu}_{\ell})\big),
\end{eqnarray}
\end{definition}
with $\theta$ the step function given by
\begin{equation}
\theta(\mu) = \left\{
\begin{array}{ll}
0 \qquad & (\mu  <   0) \\
1        & (\mu \geq 0).
\end{array}
\right.
\label{step}
\end{equation}

\subsubsection{Connection with cosets of affine algebras}
Consider the coset ${\cal C}_{n,1,1}$. The branching
functions corresponding to this coset can be defined as follows.
Let $V(\Lambda)$ be an $\widehat{sl(n)}$ HWM
with highest weight $\Lambda$, and let $|\La\rangle$ be
its highest weight vector.
Consider the tensor product decomposition
\begin{equation}
V(\Lambda_{    k}) \otimes
V(\Lambda_{i+j-k})
=
\sum_{\Lambda\in P_2^+}
\Omega_{\Lambda_k, \Lambda_{i+j-k}, \Lambda} \otimes
V(\Lambda).
\end{equation}
Among all vectors in the tensor product on the left hand
side, $\Omega_{\La_k, \La_{i+j-k},\La}$ is the
space of highest weight vectors whose weights are equal to
$\La$ mod $\Z \delta$. The connection between
$\Omega_{\La_k, \La_{i+j-k},\La}$, and $\Pset$
is as follows: It has been shown in \cite{JMMO} that
in the limit of $L \rightarrow \infty$, there is a bijection
between the set of base vectors in
$\Omega_{\La_k, \La_{i+j-k}, \La_i + \La_j}$,
and the set of paths in $\Pset$. This implies that the paths
of $\Pset$ are characterized by weights.
Under this bijection the ground-state path associated to
$\Pset$ is identified with $|\La_k\rangle \otimes
|\La_{i+j-k}\rangle \in \Omega_{\La_k,
\La_{i+j-k},\La_k+\La_{i+j-k}}$.

It turns out that the
weight of a path can be expressed in terms of its energy function
as
\begin{definition}[weight of a path $p \in \Pset$]
\begin{equation}
\hbox{\sl wt}(p) = \Lambda_i + \Lambda_j - E(p) \delta.
\label{wtp}
\end{equation}
\end{definition}

\subsubsection{Finite analogues of branching functions}
Given the above considerations, we define finite analogues
of the branching functions $B_L$ for the coset
${\cal C}_{n,1,1}$, as the generating function
of the weighted paths in $\P(\La_i+\La_j,\La_k)$,
\begin{equation}
B_L(\La_i+\La_j,\La_k) =
\sum_{p\in \P(\La_i+\La_j,\La_k)}
q^{E(p)}.
\end{equation}

\subsection{Bosonic expressions}
We are interested in expressions for the generating function
$B_L(\La_i+\La_j,\La_k)$.
In \cite{JMO}, the following {\em bosonic} expression for
$B_L(\La_i+\La_j,\La_k)$ was obtained using recurrence
relations based on the sieving method explained in \S 1:
\begin{theorem}\label{Tb}
Let $\la=\sum_{i=0}^{n-1} \la_i \,
\widehat{\imath} + \Z \delta\in P$,
with all $\la_i\geq 0$ and $\sum_{i=0}^{n-1} \la_i=N$.
For such $\la$ set
\begin{equation}
\MMult{N}{\lambda}{q} =
\frac{(q)_N}{(q)_{\la_0} \ldots (q)_{\la_{n-1}}},
\label{qmultinomial}
\end{equation}
with $(q)_m=\prod_{k=1}^m (1-q^k)$ $(m\geq 1)$ and
$(q)_0=1$.
Also, let $\cal W$ denote the Weyl group
of $\widehat{sl(n)}$ (see e.g., \cite{JMO}, p91).
Then
\begin{eqnarray}
\lefteqn{
B_L(\La_i+\La_j,\La_k) =
q^{-|\La_{i+j-k}|^2/2}} \nonumber \\
& & \nonumber \\
& & \quad \times
\sum_{w\in {\cal W}} (\mbox{det} \; w) \;
b_{L,i+j-k}(\La_k+\La_{i+j-k+L}+\rho - w(\La_i+\La_j+\rho) ),
\label{BL}
\end{eqnarray}
where
\begin{equation}
b_{L,i}(\la) = q^{|\la-\La_{i+L}|^2/2}
\MMult{L}{\la}{q}.
\label{bdef}
\end{equation}
\end{theorem}
For proof we refer the reader to \cite{JMO}.
\section{K-graphs}
Using matrices as intermediate structures, we give
an alternative representation of the weighted paths
on $P^+_2$ in terms of Ferrers graphs (or, equivalently,
Young diagrams) which satisfy certain restrictions.
We refer to these Ferrers graphs, which were introduced
and extensively studied by the Kyoto school
(see \cite{MM,JMMO} and references therein),
as K-graphs.

\subsection{Interpolating matrices}
In this subsection, we associate
a matrix $M(p)$ with 2 rows to each path $p\in \Pset$.
\begin{definition}[domain wall]\label{defwall}
Let $\iota(p)=(\mu_0,\cdots,\mu_L)$ be the integer sequence
of $p \in \Pset$.
If $\mu_{\ell}-\mu_{\ell-1}\equiv h_{\ell}+1~(0<h_{\ell}<n)$,
we say that there is a domain wall in the sequence $\iota(p)$,
of height $h_{\ell}$ at position $\ell$.
\end{definition}
Given a path $p$ with $N$ domain walls
of heights $h_1,\cdots,h_N$ at the positions
$x_1,\cdots,x_N$, respectively,  we define
the interpolating matrix $M(p)$ as
\begin{definition}[interpolating matrix]
\begin{equation}
M(p)=\left(
\begin{array}{cccc}
x_1&(x_2-x_1)&\cdots&(x_N-x_{N-1})\\
h_1&h_2    &\cdots&h_N
\end{array}
\right).
\end{equation}
\end{definition}

\medskip
\begin{example}\label{ex2}
The interpolating matrix of $p^{(2)}$ in Example~\ref{ex1} is
$$
M(p^{(2)})=\left(
\begin{array}{ccc}
1&2&3\\
3&3&2
\end{array}
\right).
$$
\end{example}

\subsection{K-graph representation of a path}
Let $p$ be a path and $M(p)$ its interpolating matrix of the
form
\begin{equation}
M(p)=\left(
\begin{array}{clcr}
w_1&w_2&\cdots&w_N\\
h_1&h_2&\cdots&h_N
\end{array}
\right).
\end{equation}
Consider a two dimensional square lattice with
an $(x,y)$-coordinate system.
Set $W=w_1+\cdots+w_N,H=h_1+\cdots+h_N$. Starting from $(0,-H)$,
we draw a polygon by moving $w_1$ steps to the right,
then $h_1$ steps up, then $w_2$ steps to the right, etc.,
until we reach the point $(W,0)$. Connecting $(0,-H)$ and $(W,0)$
with the origin by straight line-segments, the resulting graph is
{\em the Ferrers graph} or {\em Young diagram} corresponding
to the original path, see Figure~\ref{fig1sec3}(a).

\begin{definition}[K-graph]
A Ferrers graph obtained from a path $p$ on $P_2^+$, as described
above, is called a K-graph.
\end{definition}
\begin{definition}[$\G(\La_i+\La_j,\La_k)$]
$\G(\La_i+\La_j,\La_k)$ is defined as the
set of K-graphs corresponding to the set of path
$\Pset$.
\end{definition}
\begin{definition}[Profile of a graph]
The set of horizontal and vertical line segments used to
construct a K-graph form the profile of a graph.
\end{definition}

\medskip
\begin{example}
The K-graph corresponding to the interpolating matrix
of Example~\ref{ex2} is shown in Figure~\ref{fig1sec3}(b).
\end{example}

\begin{definition}[concave corner]
A corner of the form $\lceil$.
\end{definition}
\begin{definition}[convex corner]
A corner of the form $\rfloor$.
\end{definition}
\begin{definition}[plain of width $w$]
A horizontal line segment of $w$ nodes (or boxes) marked by a
concave corner to its left and convex corner to its right.
\end{definition}
\begin{definition}[cliff of height $h$]
A vertical line segment of $h$ nodes (or box\-es) marked by a
convex corner at its bottom and a concave corner at its top.
\end{definition}
Notice that a cliff on a K-graph corresponds to a domain wall
in the corresponding integer sequence.

\begin{remark}
{}From now on, we concentrate on K-graphs
in $\G(2\La_0,\La_k)$,
unless otherwise stated.
\end{remark}

\subsubsection{From a graph to its sequence of integers}
For a graph $G\in \Gset$, we can recover the corresponding
integer sequence $\iota(p)$ as follows.
Let
\begin{equation}
M=\left(
\begin{array}{clcr}
w_1&w_2&\cdots&w_N\\
h_1&h_2&\cdots&h_N
\end{array}
\right)
\label{Mwh}
\end{equation}
be the interpolating matrix corresponding to $G$. Set
$H=h_1+\cdots h_N$, and take the integer sequence
$(0,1,2,\ldots,n-1,0,1,\ldots,n-1,0,1,\ldots)$
of length $H+L+1$.
Strarting from the left moving to the right,
we now keep the first $w_1$ integers,
then remove the next $h_1$ integers, then keep the next $w_2$
integers, remove the next $h_2$ integers, etc.
The remaining sequence of exactly $L+1$ integers corresponds
to $\iota(p)$.

\subsubsection{From a graph to a path}
To go from a graph $G\in \Gset$,
to its corresponding path $p=(\la_0,\ldots,\la_L)$
on $P_2^+$
we simply first construct the sequence of integers
$\iota(p)=(\mu_0,\ldots,\mu_L)$ as described above.
We then compute $\la_{\ell+1}=\la_{\ell}+\widehat{\mu}_{\ell}$
using $\la_0=2\La_0$.

\subsubsection{Conditions on $\Gset$}\label{seccond}
Among all K-graphs, those in $\Gset$ are characterized
by the following conditions:
\begin{description}
\item[K1] $W \leq L$, with $W$ the number of nodes in the first row.
\item[K2] $H+k \equiv 0$, with $H$ the number of nodes in the
          first column.
\item[K3] $h_{i-1}+w_i+h_i\equiv 0$ and $0<h_i<n$
          for $1\le i\le N$,
          with $h_0=n$.
\end{description}
{\bf K1} is obvious. {\bf K2} is obtained by considering the
$L$-th component of $\iota(p)$ and the boundary condition.
To obtain {\bf K3}, suppose the $i$-th cliff occurs
at the $r$-th position. We can assume $\la_r=\La_a+\La_b$,
$\la_{r-w_i}=\La_a+\La_{b-w_i}$ for some $a,b$. Now we have
$\iota(p)_{r-1}\equiv b-1$, $\iota(p)_{r-w_i}\equiv b-w_i$.
Since there are cliffs at the $r$-th and $(r-w_i)$-th position,
we should have $\iota(p)_r\equiv a$ and $\iota(p)_{r-w_i-1}\equiv a-1$.
Thus we get $h_i\equiv a-b$, $h_{i-1}\equiv b-w_i-a$, which
gives {\bf K3}.

\begin{figure}[hbt]
\centerline{\epsffile{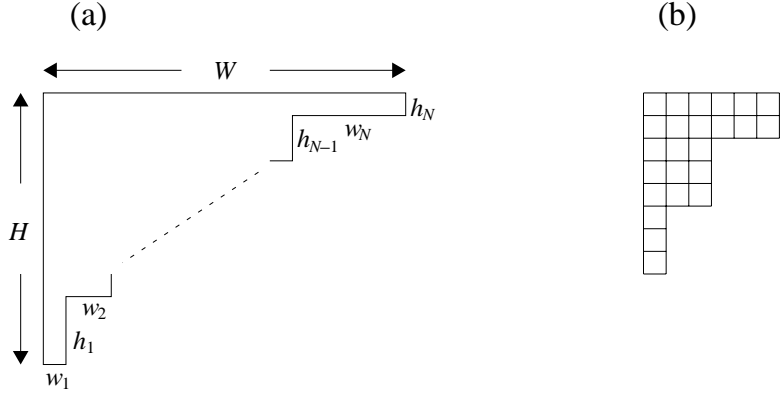}}
\caption{(a) The general form of a K-graph.
         (b) A K-graph in ${\cal G}_6(2\La_0,\La_0)$ for $n=4$.}
\label{fig1sec3}
\end{figure}

\subsection{Fermionic expressions}
We now wish to calculate the following sum:
\begin{equation}
F_L(\La_i+\La_j,\La_k)=\sum_{G\in\G(\La_i+\La_j,\La_k)}q^{|G|/n},
\label{fer}
\end{equation}
where $|G|$ denotes the number of nodes in $G$.
Regarding the above, we have the following theorem \cite{MM,JMMO}.
\begin{theorem}
Let $p$ be a path in $\Pset$, and $G(p)$ the corresponding
K-graph. The number of nodes of $G(p)$ is given by
\begin{equation}
|G|=\sum_{\ell=0}^{n-1}m_{\ell},
\end{equation}
where $m_{\ell}$ is determined from
\begin{equation}
(\La_k+\La_{i+j-k})-\hbox{\sl wt}(p)=\sum_{\ell=0}^{n-1}
m_{\ell}\alpha_{\ell}.
\end{equation}
\end{theorem}

Using
$$
( \La_i\mid \alpha_j ) =\delta_{ij}\ \  (i,j=0,...,n-1),
$$
we obtain
$$
m_l=(\La_l \mid (\La_k+\La_{i+j-k})-wt(p))\ \ (l=0,...,n-1).
$$

Since we define the sum (\ref{fer}) in the ``principal picture'', \ie each
node has equal weight $1/n$, it is invariant under the Dynkin
diagram automorphisms. Thus we can reduce the calculation of (\ref{fer})
to that of $F_L(\La_0+\La_j,\La_k)$. From now on, we hence assume $i=0$.

Setting $(\La_k+\La_{j-k})-(\La_0+\La_j)=
\sum_{\ell=0}^{n-1}\bar{m}_{\ell}\alpha_{\ell}$,
we have $\sum_{\ell=0}^{n-1}\bar{m}_{\ell}=k(j-k)$
for $j\geq k$, $=(k-j)(n-k)$
for $j<k$. Calculating $|\La_k|^2+|\La_{j-k}|^2-|\La_j|^2$ and comparing
(\ref{fer}) with the bosonic expression, we obtain
\begin{equation}
F_L(\La_0+\La_j,\La_k)=q^{(|\La_k|^2+|\La_{j-k}|^2-|\La_j|^2)/2}
B_L(\La_0+\La_j,\La_k).
\label{BtoF}
\end{equation}

In the remainder of this paper we will compute a fermionic type
of expression for $F_L$. Given (\ref{BtoF}) and the bosonic
expression (\ref{BL}) for $B_L$, this gives rise to polynomial
identities for the finite analogues of the branching functions
of the coset ${\cal C}_{n,1,1}$.
\section{Parents}
{}From the set of all K-graphs in $\Gset$
we select a subset of graphs to be called parent graphs,
or simply parents.
Let $\vec{m}^t=(m_1,\cdots,m_{n-1})\in(\Znn)^{n-1}$,  such that
\begin{equation}
k+\sum_{i=1}^{n-1} i m_i \equiv 0
\label{column1}
\end{equation}
and let $M$ be the interpolating matrix of a graph $G\in \Gset$,
with entries
\begin{equation}
M(p)=\left(
\begin{array}{clcr}
w_1&w_2&\cdots&w_N\\
h_1&h_2&\cdots&h_N
\end{array}
\right).
\end{equation}
\begin{definition}[parent associated to $\vec{m}$]
$G$ is the parent associated to $\vec{m}$ if
\begin{equation}
\left\{
\begin{array}{rcl}
h_1,\ldots,h_{m_{n-1}}&=&n-1 \\
h_{m_{n-1}+1},\ldots,h_{m_{n-1}+m_{n-2}}&=& n-2 \\
& \vdots & \\
h_{N-m_1+1}, \ldots, h_N &=& 1,
\end{array} \right.
\end{equation}
with $N=\sum_{i=1}^{n-1} m_i$, and
\begin{equation}
h_{i-1}+ w_i + h_i = 2n \qquad 1\leq i \leq N,
\label{hwh2n}
\end{equation}
where we recall that $h_0 = n$.
\end{definition}

\begin{example}
The K-graph of $p^{(2)}$ shown in Figure~\ref{fig1sec3}(b)
is the parent associated to $\vec{m}^t=(0,1,2)$.
\end{example}

\subsubsection{The number of nodes of a parent graph}
Let $C$ be the Cartan matrix of $sl(n)$, \ie
$C_{ij}=2\delta_{i,j}-\delta_{i,j-1}-\delta_{i,j+1}\,(i,j=1,\cdots,n-1)$.
The inverse of $C$ is then given by the following formula:
\begin{equation}
(C^{-1})_{ij}=
\left\{
\begin{array}{ll}
\frac{i(n-j)}{n}&(i\le j)\\
& \\
\frac{j(n-i)}{n}&(i>j).
\end{array}
\right.
\label{invC}
\end{equation}
With this definition we have the following lemma:
\begin{lemma}
The number of nodes of the parent associated to $\vec{m}$
is given by
$n\, \vec{m}^tC^{-1}\vec{m}$.
\end{lemma}
Though the proof of this statement is rather elementary, we need
to take some care as some of the entries of $\vec{m}$ can actually
be zero. In the following we use the notation $<i>$ to denote
$\sum_{j=1}^i m_{n-j}$. Clearly, $<i>-<i-1>=m_{n-i}$.
We now compute the number of nodes of a parent $N(\vec{m})$ as follows
\begin{equation}
N(\vec{m})= \sum_{i=1}^N \sum_{j=1}^i h_i w_j
= \sum_{k=1}^{n-1} \sum_{i=1+<k-1>}^{<k>}
\left(
\sum_{\ell=1}^{k-1} \sum_{j=1+<\ell-1>}^{<\ell>} + \sum_{j=1+<k-1>}^i
\right) h_i w_j.
\end{equation}
Now use (\ref{hwh2n}) and $h_{1+<i-1>}=\ldots= h_{<i>} = n-i$ to get
\begin{eqnarray}
N(\vec{m})&=&
\sum_{k=1}^{n-1} \sum_{i=1+<k-1>}^{<k>} (n-k) \nonumber \\
& & \qquad \times
\left(
\sum_{\ell=1}^{k-1} \sum_{j=1+<\ell-1>}^{<\ell>}\left\{
2 \ell -(h_{<\ell-1>}-(n-\ell))\delta_{j,<\ell-1>+1} \right\}
\right. \nonumber \\
& & \quad \qquad  \left.
 + \sum_{j=1+<k-1>}^i
\left\{
2 k -(h_{<k-1>}-(n-k))\delta_{j,<k-1>+1} \right\} \right).
\end{eqnarray}
Finally, after some changes of variables, we obtain
\begin{eqnarray}
N(\vec{m})&=& 2 \sum_{k=1}^{n-1} \sum_{\ell=1}^{k-1}
\ell (n-k)  m_k m_{\ell} \nonumber \\
& & + \sum_{k=1}^{n-1} (n-k) \sum_{i=1}^{m_{n-k}}
\left( 2 i k
- \sum_{\ell=1}^k \sum_{j=1}^{m_{n-\ell}}
( h_{<\ell-1>}-(n-\ell)) \delta_{j,1} \right) \\
&=&
2 \sum_{k=1}^{n-1} \sum_{\ell=1}^{k-1}
\ell (n-k) m_k m_{\ell}
+ \sum_{k=1}^{n-1} \sum_{i=1}^{m_{n-k}} k(n-k)(2i-1).
\nonumber
\end{eqnarray}
Summing over $j$ and $i$ this results
in $n\,\vec{m}^t C^{-1} \vec{m}$.
\section{g-components}
Now that we have distinguished a subset of all
K-graphs as parents, we wish to describe the
minimal connected configuration of nodes that
can be removed or added to a K-graph in
$\G(2\La_0,\La_k)$ to obtain another K-graph in
$\G(2\La_0,\La_k)$. Since, as we will see in \S 7,
these configurations are generated by Gaussian
polynomials, we call them g-components. Eventually,
we will show that those graphs which are related by
addition and removal of g-components belong to the
same {\em sector}, and we will use this observation
to relate any non-parent graph to a parent graph.

\begin{figure}[hbt]
\centerline{\epsffile{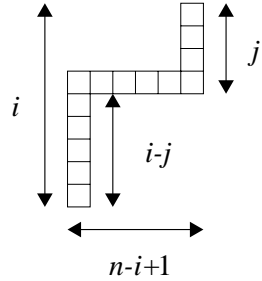}}
\caption{An $(i, j)$-component}
\label{gcomp}
\end{figure}

\begin{definition}[$(i,j)$-component]
For all $i=1,\ldots,n-1$ and all $j=1,\ldots,i$,
we define an $(i,j)$-component as a connected
configuration of $n$ nodes, as shown in Figure~\ref{gcomp}.
\end{definition}
Some important characteristics of an $(i,j)$-component are:
\begin{enumerate}
\item[G1] It consists of $n$ nodes.
\item[G2] It has total height $i$.
\item[G3] It has total width $n-i+1$.
\item[G4] it has (at most) two cliffs, one (the lower) of height
          $i-j$, and one (the upper) of height $j$.
\end{enumerate}
We further note that for an $(i,i)$-component the lower
cliff vanishes resulting in a configuration with a single cliff.

\begin{definition}[$i$-component]
An $(i,j)$-component for arbitrary $j$ is called an $i$-component.
\end{definition}

\begin{definition}[g-component]
An $(i,j)$-component for arbitrary $i$ and $j$ is called a
g-component.
\end{definition}

We are now interested in the addition/removal of
g-components to/from a K-graph.
Clearly, in adding or removing a g-component  to or from a K-graph
in $\G(2\La_0,\La_k)$, we demand that the resulting
graph is again a graph in $\G(2\La_0,\La_k)$.
However on top of this we impose one additional condition, which
basically defines our sectors.

\subsection{Removing an $i$-component}\label{Raic}
\begin{figure}[hbt]
\centerline{\epsffile{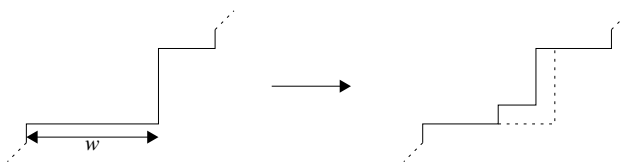}}
\caption{Removing an $(i,i)$-component is only allowed when
$w>2(n-i)$.  The extra dotted lines in the
resulting graph are to indicate the nodes which are removed.}
\label{removal}
\end{figure}

The removal of an $(i,j)$-component from a K-graph in
$\G(2\La_0,\La_k)$ is allowed provided the following two
conditions are satisfied:
\begin{enumerate}
\item[R1]
The resulting graph is again a K-graph in
$\G(2\La_0,\La_k)$.
\item[R2]
If $j=i$ as in Figure~\ref{removal}, we demand that $w>2(n-i)$.
\end{enumerate}

\begin{definition}[$i$-candidate]
An $i$-component one is allowed to remove from a K-graph is
called an $i$-candidate:
\end{definition}
Since for any K-graph in $\G(2\La_0,\La_k)$ we have
$h_{j-1}+w_j+h_j \equiv 0$, three kinds of candidates can occur.
\begin{enumerate}
\item
$h_{j-1}+w_j+h_j=n$ and $w_{j-1}>1$.
In this case we can remove an $(h_{j-1}+h_j,h_j)$-component.

\item
$h_{j-1}+w_j+h_j=2n$ and $w_j>2(n-h_j)$.
In this case we can remove
an $(h_j,h_j)$-component.

\item
$h_{j-1}+w_j+h_j\geq 3n$.
In this case we can remove an $(h_j,h_j)$-component.

\end{enumerate}

Scanning the profile of a non-parent graph, several $i$-candidates
may occur.
\begin{definition}[leading $i$-candidate]
The leading $i$-candidate is the down- and left-most $i$-candidate.
\end{definition}

\subsection{Attaching an $i$-component to a graph}
Attaching an $i$-component to a K-graph in
$\G(2\La_0,\La_k)$ is allowed provided the following
conditions are satisfied:
\begin{enumerate}
\item[A1]
The resulting graph is again a K-graph in
$\G(2\La_0,\La_k)$.
\item[A2]
We do not generate an $i'$-candidate, with $i'>i$.
\end{enumerate}

\begin{definition}[$i$-vacancy]
An $i$-vacancy is a position on the profile such that one is
allowed to attach
an $i$-component.
\end{definition}

\begin{figure}[hbt]
\centerline{\epsffile{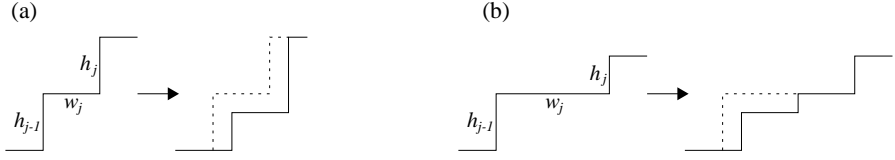}}
\caption{(a) Attaching an $(h_{j-1}+h_j,h_j+1)$-component with $h_j \neq 1$.
(b) Attaching an $(h_{j-1},1)$-component.
The extra dotted lines in the
resulting graphs are to indicate the profile before attaching the
g-component.}
\label{attach}
\end{figure}

An important statement about $i$-vacancies is the following.
Given a sequence $\{\mbox{cliff,plain,cliff}\}$ of dimensions
$h_{j-1},w_j,h_j$ such that $h_{j-1}+w_j+h_j=n$ and
$h_{j-1}+h_j=i$.
Then the following holds:
\begin{lemma}\label{annil}
If $w_{j+1}=n-i$, the above sequence
is {\em} not an $i$-vacancy.
\end{lemma}
To proof this assume the above sequence is an $i$-vacancy.
Hence we can attach an $(i,h_j+1)$-component as shown in
Figure~\ref{attach}(a).
Note that in doing so the height of the $j$-th cliff increases
by one to $h_j'=h_j+1$, and the width of the $(j+1)$-th plain
decreases by one to $w_{j+1}'=w_{j+1}-1=n-i-1$.
Thus we compute
\begin{equation}
0<h'_j+w'_{j+1}+h_{j+1}=h_j+n-i+h_{j+1}\leq i+(n-i)+(n-1) =2n-1.
\end{equation}
Since $h'_j+w'_{j+1}+h_{j+1}\equiv n$, we conclude that
$h'_j+w'_{j+1}+h_{j+1}=n$, and $h'_j+h_{j+1}=i+1$.
But these are the characteristics of an $(i+1)$-candidate.
By the second condition for attaching $i$-components
this contradicts our assumption that the initial sequence was
an $i$-vacancy.

With the above lemma we note that two kinds of vacancies
may occur. The first occurs if we have a sequence
$\{\mbox{cliff,plain,cliff}\}$
of dimensions
$h_{j-1},w_j,h_j$
such that
$h_{j-1}+w_j+h_j=n$
and
$w_{j+1}>n-h_j-h_{j+1}$.
In this case we can always attach a $(h_{j-1}+h_j,h_j+1)$-component,
as shown in Figure~\ref{attach}(a).
The second occurs if we have have a sequence
$\{\mbox{cliff,plain,cliff}\}$
of dimensions
$h_{j-1},w_j,h_j$
such that
$h_{j-1}+w_j+h_j\geq 2n$.
In this case we can always attach an $(h_{j-1},1)$-component,
see Figure~\ref{attach}(b).

\vspace{5mm}

\section{Descendants}
In previous sections, we classified all admissible
K-graphs into parents and non-parents. We need to show
that each non-parent is a {\em descendant} of a unique
parent.  More precisely, we show that:
\begin{enumerate}
\item Given a non-parent graph, there is a reduction
      procedure, such that one can reduce it to a unique
      parent graph.
\item The reduction procedure is reversible: given a parent
      graph, there
      is a composition procedure to recover the original
      non-parent graph.
\end{enumerate}

Because the reduction procedure is reversible, any non-parent
graph is a descendant of a unique parent graph. Thus the
set of all admissible K-graphs can be divided into
non-overlapping sectors. Each sector contains, and is labelled
by a parent graph. Any admissible K-graph belongs to one
and only one sector.

\subsection{Reducing non-parent graphs}
Given a non-parent K-graph, we can reduce it to a parent
graph as follows.
\begin{enumerate}

\item[Red0] Set $i=n-1$.

\item[Red1] Search for the leading $i$-candidate
          and, if it exists, remove it.

\item[Red2] Repeat the above step till no more
          $i$-candidates are found.

\item[Red3] Set $i\to i-1$ and, if $i\geq 1$, repeat Red1-Red3.
\end{enumerate}

To prove that a reduced
graph is indeed a parent, we proceed as follows:
Consider a profile with a sequence $\{\mbox{cliff,plain,cliff}\}$
of dimensions $h_{j-1},w_j,h_j$, respectively.
Suppose that the part of the profile below the
above sequence belongs to a parent, i.e.,
$h_{k-1}+w_k+h_k=2n$ and $h_{k-1}\geq h_k$
for $k=1,\ldots,j-1$. We wish to show that if the
above sequence does not represent a candidate, it belongs to a parent.
{}From section~\ref{Raic} we see that unless
$h_{j-1}+w_j+h_j=2n$ and $w_j\leq 2(n-h_j)$
or
$h_{j-1}+w_j+h_j=n$ and $w_{j-1}=1$,
we always have a candidate.

In the first case we get
$2n=h_{j-1}+w_j+h_j\leq h_{j-1}-h_j + 2n$ and thus
$h_{j-1}\geq h_j$. This is
precisely the right sequence for a parent and we get
$h_{k-1}+w_k+h_k=2n$ and $h_{k-1}\geq h_k$
for $k=1,\ldots,j$.
The second case can in fact never occur.
Since $h_{j-2}+w_{j-1}+h_{j-1}=2n$ and $h_{j-2},
h_{j-1}<n$ we find that $w_{j-1}>1$.

\subsection{Generating descendants from parents}
Given the parent associated to $\vec{m}$.
Each cliff of height $i$
plus the plain immediately to the right
of this cliff, forms an $i$-vacancy.
Hence we have $m_i$ $i$-vacancies.
To obtain an arbitrary descendants of the parent under consideration,
we proceed as follows.
\begin{enumerate}
\item[Gen0] Set $i=1$.

\item[Gen1] Set $j=1$.

\item[Gen2] Attach $k_j^{(i)}$ $i$-components
          to the $j$-th $i$-vacancy
          counted from the right.

\item[Gen3] Set $j \to j+1$.
          If $j\leq m_i$, go to Gen2.
          If $j=m_j+1$, set $i\to i+1$, and if
          $i\leq n-1$, go to Gen1.
\end{enumerate}

To properly interpret these rules some important
remarks need to be made.
First, when we say ``attach $k_j^{(i)}$  $i$-components
to the $j$-th $i$-vacancy'' this should be understood as follows.
Attaching an $i$-component to an $i$-vacancy has the effect
of moving the vacancy to the right. Hence attaching the $k$-th
$i$-component means attaching an $i$-component to the image
of the $i$-vacancy after attaching the $(k-1)$-th $i$-component.
Second, it may occur that attaching an $i$-component to
an $i$-vacancy does not have the effect of moving the
$i$-vacancy to the right, but annihilates the vacancy.
Hence, there are bounds on the numbers $k_j^{(i)}$.
In the next section we will show that these bounds are
as follows:
\begin{equation}
0\leq k_{m_i}^{(i)}\leq \ldots \leq k_2^{(i)} \leq k_1^{(i)} \leq \ell_i,
\end{equation}
with $\ell_i$ fixed by (\ref{ldef}).

\subsection{Reversibility}
Remains the proof that our rules for attaching and removing
g-components are reversible.
This is true by construction.

\section{Proof of Gaussians}
In this section we prove that for the case of
${\cal G}_L(2\La_0,\La_k)$, the generating function
for attaching the $i$-components to the
parent graph associated to $\vec{m}$,
is given by the Gaussian
polynomial
\begin{equation}
\Mult{\ell_i+m_i}{m_i}{q},
\end{equation}
where
\begin{equation}
\Mult{N}{m}{q} = \left\{
\begin{array}{ll}
\displaystyle \frac{(q)_N}{(q)_m (q)_{N-m}} \qquad & 0\leq m\leq N \\
& \\
0 & \mbox{otherwise.}
\end{array}
\right.
\end{equation}
and
\begin{equation}
\vec{\ell} = C^{-1} \left(
L \, \vec{\e}_{n-1} + \vec{\e}_r -2 \vec{m} \right),
\label{ldef}
\end{equation}
and $\vec{\e}_i$ the $(n-1)$-dimensional unit vector with entries
$\left(\vec{\e}_i\right)_j = \delta_{i,j}$ and
with $0< r\leq n$ fixed by
\begin{equation}
L-2k\equiv r.
\label{rdef}
\end{equation}

\subsection{The $r=n$ case}
To prove the above result we first treat the simpler
case of $r=n$. In the next subsection we then show how
to modify this to obtain (\ref{ldef}) for general $r$.

We start with the following important fact, used extensively throughout
this section:
\begin{lemma}\label{lem7}
For $G\in{\cal G}_L(2\La_0,\La_k)$,
let $W$ be the number of nodes in the first row
and $h_N$ the height of the $N$-th (uppermost) cliff,
see Figure~\ref{fig1sec3}(a).
Then
\begin{equation}
W-h_N\equiv 2k.
\label{Wmh}
\end{equation}
\end{lemma}
We proof this by implementing the conditions K2 and K3
of section~\ref{seccond},
defining the K-graphs in $\Gset$.
Recalling that $W$ is the number of nodes in the
first row of a K-graph, we have
\begin{eqnarray}
W&=&\sum_{i=1}^N w_i \equiv -\sum_{i=1}^N (h_{i-1}+h_i)
\equiv h_N - 2 \sum_{i=1}^N h_i  \nonumber \\
&=& h_N - 2 H \equiv h_N+2 k,
\end{eqnarray}
which proves our claim.

\subsubsection{$i$-strips}\label{sgs}
\begin{figure}[bt]
\centerline{\epsffile{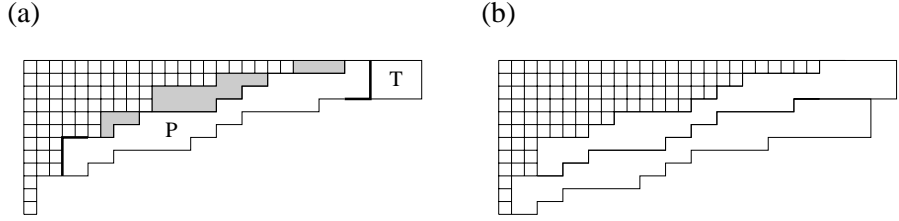}}
\caption{(a) The construction of the first $3$-strip
for a typical K-graph $G$. The part of $G$ corresponding to the
parent graph is drawn with open nodes/boxes and the
(already) placed $i'$-components $(i'=1,2)$ are shown in grey.
The bold segments are the left- right terminals of
the principal region (marked P). The tail of the strip is marked with T.
(b) The two $3$-strips for the K-graph in (a).}
\label{figgstrip}
\end{figure}

We are interested in the placement of the $i$-components.
{}From the rules for placing the latter, it is natural to
define $m_i$ {\em $i$-strips} as follows:

\begin{enumerate}
\item[S1]

We define the $i$-strip as consisting of two regions:
a principal region, and a tail. The principal region is defined
in terms of a top segment, a bottom segment, and a left, and
a right segment. The left and right segments will be called
left and right terminal. We start by defining the principal
region.

Consider the profile $P_i$ of the K-graph after attaching
all $i'$-components, with $i'=1, 2, \cdots, i-1$,
but before attaching any component of height $i$ or
higher.

If the $N$-th (highest) cliff of $P_i$ has height $h_N$,
extend the {\it ceiling} of $P_i$ by drawing a horizontal
line of width $i-h_N$ starting from the top-right corner
of the right-most node of the top row of the graph, and
extending to the right. This line will serve as a ceiling
to the $i$-components that we will add shortly.

Consider the segment of $P_i$ $2(n-i)$ columns to the right of the
right-most cliff of height $i$.
This will be the top
segment of the $i$-strip. Let us denote this segment by $P_i^0$.

Now we proceed to define the bottom segment of the $i$-strip.
Move $P_i^0$ to the left by $2(n-i)$ columns, and downwards by
$i$ rows. Denote this shifted profile by $P_i^1$. This is
the bottom segment that we are looking for.

Finally, close the figure formed by the top and bottom segments
as follows: draw a plain of width $2(n-i)$ followed by a cliff
of height $i$ to the right (left) of $P_i^0 (/P_i^1)$, called
the right (left) terminal, respectively. As a result, we now
have a region enclosed by $P_i^0$, $P_i^1$ and the left and right
terminals. This defines the {\em principal region} of the first
$i$-strip.

Next, we define the tail of the first $i$-strip as follows:
Compute $M \equiv 0$ to be the largest integer such that
$W+i-h_N+M\leq L$. The tail of the $i$-strip is a rectangle,
of width $M$ and height $i$, that we place to the right of the
principal region in the first row.

The principal region plus the tail define the complete first
$i$-strip. An example of the first $3$-strip in a typical
K-graph for $\widehat{sl(4)}$ is shown in Figure~\ref{figgstrip}(a).

\item[S2]
We draw the second $i$-strip by simply shifting
the first $i$-strip to the left and down by
$2(n-i)$ columns and $i$ rows, respectively.

\item[S3]
We repeat the step S2 $(m_i-1)$ times.
That is, we define the $(j+1)$-th $i$-strip
by translating the $j$-th $i$-strip
to the left by $2(n-i)$ columns and downwards
by $i$ rows.
In Figure~\ref{figgstrip}(b) we have shown the construction
of the $3$-strips for a typical
example of a K-graph
for $\widehat{sl(4)}$.
\end{enumerate}

By construction, adding the $i$-components corresponding to
the $j$-th $i$-va\-can\-cy (counted from the right),
corresponds to filling the
$j$-th $i$-strip from left to right.
In constructing an arbitrary descendant, we will not necessarily
fill the complete $j$-th $i$-strip. Furthermore, we will show
below that the filling of the $(j+1)$-th strip is bound by
the degree of filling of the $j$-th strip. In particular, we will
show that if the $j$-th strip is filled with $k_j^{(i)}$ $i$-components,
then the $(j+1)$-th strip cannot be filled with more than $k_j^{(i)}$
$i$-components. Since each $i$-component contains $n$ nodes, and
thus contributes a single factor $q$, we obtain the following
expression for the generating function attaching the $i$-components:
\begin{equation}
G_i(q) = \sum_{k_1^{(i)}=0}^{\ell_i} \sum_{k_2^{(i)}=0}^{k_1^{(i)}} \ldots
\sum_{k_{m_i}^{(i)}=0}^{k_{m_i-1}^{(i)}}
q^{k_1^{(i)}+k_2^{(i)}+\ldots +k_{m_i}^{(i)}}.
\end{equation}
Here the number $n \ell_i$ is
the area (=number of nodes) of the first $i$-strip.

As defined above, $G_i$ can be interpreted as the generating
function of all partitions with largest part $\leq \ell_i$ and
number of parts $\leq m_i$. Therefore
\begin{equation}
G_i(q) = \Mult{\ell_i+m_i}{m_i}{q}.
\end{equation}
Before ending this subsection, let us return to Lemma~\ref{lem7}.
We have stated above that attaching $i$-components
corresponding to the right-most $i$-vacancy corresponds
to  filling the first $i$-strip.
However, some caution needs to be taken, since in constructing the
principal region of the
first $i$-strip we have extended the profile of the
K-graph by drawing a plain of width $i-h_N$
in the first row to the right of the $N$-th cliff.
This clearly can only be done for all $i=1,\ldots,n-1$,
if $L-W\geq n-1-h_N$.
If $L_s$ is the smallest possible value of
$L$ for which a K-graph of width $W$ is possible, \ie
\begin{equation}
L_s=W+x \qquad 0\leq x \leq n-1
\end{equation}
with $x$ fixed by (\ref{rdef}),
we have
\begin{equation}
x=L_s - W = (L_s-2k)-(W-2k) \equiv r-h_N,
\end{equation}
where we have used lemma~\ref{lem7}
and the definition (\ref{rdef}) of $r$.
Since we require $x$ to be at least $n-1-h_N$
we should thus have that $r=n$ or $n-1$.
For simplicity we now assume $r=n$.

\subsubsection{Calculation of $\ell_i$}
To calculate the area of the first $i$-strip, we use the simple
property that the area remains unchanged by deforming the
strip by removing nodes from below and adding them from above.

\begin{figure}[hbt]
\centerline{\epsffile{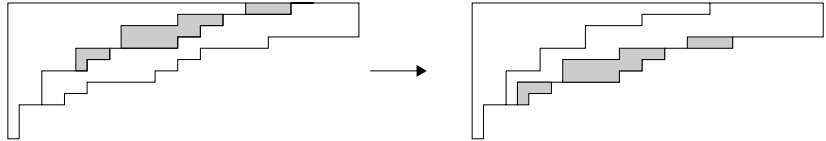}}
\caption{The deformation of the first $3$-strip of the
K-graph of Figure~6,
yielding the first $3$-strip of its parent graph.}
\label{figdeform}
\end{figure}

We now choose to deform the $i$-strip such that its
upper-side corresponds to the profile of its parent graph,
being labelled by $\vec{m}$.
For the example of Figure~\ref{figgstrip} this is shown
in Figure~\ref{figdeform}.

{}From this particular choice of deformation we can simply
compute the area over $n$ as
\begin{equation}
\ell_i(r=n)=2 \sum_{j=1}^{i-1} j m_{i-j}  +  i
\lfloor
( L + n-i - 2 \sum_{j=1}^{n-1} j m_{n-j} )/n
\rfloor,
\label{ell0}
\end{equation}
with $\lfloor x \rfloor$ denoting the integer part of $x$.
Here the first term corresponds to the area of the principal
region of the {\em deformed} $i$-strip
and the second term to the area of the tail of the deformed strip.
In particular, to compute the former we use
the fact that
it takes $(j+1)$ $i$-components to move an
$i$-vacancy (of the type shown in Figure~\ref{attach}(b))
upwards across a plain of width
$2n-2i+j$, $(j\geq 1)$.  To compute the latter,
we compute $L-(W+i-h_N)$,
using the result (\ref{Wmh}).

Recasting the definition (\ref{invC}) of the inverse Cartan matrix as
\begin{equation}
( C^{-1} )_{i,j} = \frac{i(n-j)}{n} - \sum_{p=1}^{i-1} (i-p) \delta_{j,p}
\end{equation}
and using the mod properties (\ref{column1}) and (\ref{rdef})
we thus get,
\begin{eqnarray}
\ell_i(r=n) &=&  \frac{i L}{n} -2 \left( \sum_{j=1}^{n-1} \frac{i(n-j)}{n}-
\sum_{j=1}^{i-1} (i-j) \right) m_j \nonumber \\
& & \nonumber \\
&=&
L(C^{-1})_{n-1,i}-2(C^{-1} \vec{m})_i .
\end{eqnarray}
This proves (\ref{ldef}) for $r=n$.

\vspace{5mm}

\subsubsection{Proof of Gaussian form}
Remains the proof that the filling of the
$(j+1)$-th $i$-strip is bound by that
of the $j$-th $i$-strip.

\begin{figure}[hbt]
\centerline{\epsffile{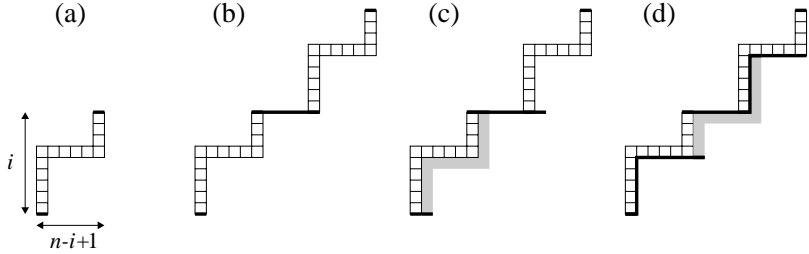}}
\caption{(a) A typical $i$-component.
The bold lines indicate part of the boundary of the
$i$-strip. (b) The $k_j^{(i)}$-th $i$-component
in both the $j$-th and the $(j+1)$-th $i$-strip.
(c) and (d) The two forbidden placements of an additional
$(k_j^{(i)}+1)$-th $i$-component (shown in grey)
in the $(j+1)$-th strip.}
\label{fig5}
\end{figure}

Let us assume that we have placed $k_j^{(i)}$ $i$-components
in the $j$-th $i$-strip, and that the $k_j^{(i)}$-th
component is an $i$-component of the form
depicted in Figure~\ref{fig5}(a).
Let us further assume that we have already
filled the $(j+1)$-th strip with $k_j^{(i)}$ $i$-components
(this is of course always possible).
Since the $(j+1)$-th strip has identical shape as the
$j$-th strip, but is translated to the left and down
by $2(n-i)$ columns and $i$ rows,
we have the configuration shown in Figure~\ref{fig5}(b).

Our claim is now that upon attaching the $k_j^{(i)}$-th
$i$-component  in the $(j+1)$-th strip we have annihilated
the corresponding $i$-vacancy.
To see this we consider two cases.
Either the boundary separating the strips extends at least one
more entry to the right, see Figure~\ref{fig5}(c),
or the boundary progresses upwards as in~\ref{fig5}(d).
In the first case the lowest sequence $\{\mbox{cliff,plain,cliff}\}$
could be an $i$-vacancy, but since the plain
immediately above has width $n-i$, this is not the case
thanks to Lemma~\ref{annil}. Hence placing the
$i$-component in \ref{fig5}(c) as shown in grey is not allowed.
In the second case,
the middle sequence $\{\mbox{cliff,plain,cliff}\}$
could be an $i$-vacancy, but again the plain immediately
above has width $n-i$ and we can once more apply
Lemma~\ref{annil}. Hence also the placement
as shown in \ref{fig5}(d) is forbidden.

\subsection{The general $r$ case}
\begin{figure}[hbt]
\centerline{\epsffile{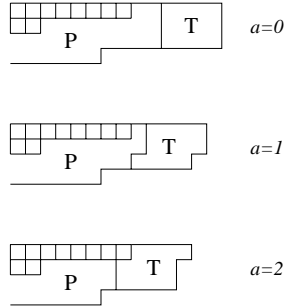}}
\caption{Reduction of the $3$-strips of Figure~6 .
Recalling that $i=r+a$, we need for $i=3$ and $n=4$
to consider the cases
$a=0,1,2$.}
\label{figred}
\end{figure}

As remarked at the end of section~\ref{sgs},
only for $r=n$ and $n-1$ we can always draw a plain of width
$i-h_N$ to the right of the $N$-th cliff without
violating the condition $W+i-h_N\leq L$, for any $i$.
If $L-2k\equiv r$, we can still do so for all $i\leq r$.
Hence for these cases the principal region
of the $i$-strips can still be defined as in section~\ref{sgs}.
However for
$i=r+a$ $(a>0)$, we have to reduce the
principal region P by removing the part of P which would
be occupied by the last $a$ components to be attached,
if P were to be completely filled
from left to right. Of, course, in this case the tail
no longer is a rectangle, but has a profile of 2 plains
and 2 cliffs.
An example of this reduction is shown in
Figure~\ref{figred}.

The above considerations
lead to the following simple modification of (\ref{ell0}):
\begin{eqnarray}
\ell_i &=&
2 \sum_{j=1}^{i-1} j m_{i-j}
- \sum_{p=1}^{n-r-1} (n-r-p) \delta_{n-i,p} \nonumber \\
& & +  i
\lfloor (
L + n-i +\sum_{p=1}^{n-r-1}(n-r-p) \delta_{n-i,p}
- 2 \sum_{j=1}^{n-1} j m_{n-j} )/n \rfloor \nonumber \\
&=&
2 \sum_{j=1}^{i-1} j m_{i-j}
- \sum_{p=1}^{n-r-1} (n-r-p) \delta_{n-i,p} +
\frac{i}{n}
\left(
L + n-r - 2 \sum_{j=1}^{n-1} j m_{n-j} \right) \nonumber \\
&=& \ell_i(r=n) + \left( \frac{i(n-r)}{n} -
\sum_{p=1}^{n-r-1}(n-r-p) \delta_{n-i,p} \right)
\label{allr} \\
&=& \ell_i(r=n) +
(C^{-1})_{r,i}, \nonumber
\end{eqnarray}
which proves the claim (\ref{ldef}).
Here we note that to obtain the first line of (\ref{allr})
one not only has to subtract the term $\sum_p (n-r-p) \delta_{n-i,p}$
to account for the reduction of the principal region, but also
to add this same term within the $\lfloor . \rfloor$. This
occurs since the
effective length available for the tail of the $i$-strips has of course
increased by the decrease of the principal region, see
Figure~\ref{figred}.

\vspace{5mm}

\subsection{Fermionic form for $F_L(2\La_0,\La_k)$}
We now have computed the number of nodes of the parent
associated to $\vec{m}$ as well as the generating function
for adding the g-components to this parent.
Collecting these two results, we obtain the following expression for
the generating function $F_L(2\La_0,\La_k)$ of K-graphs in
${\cal G}(2\La_0,\La_k)$:
\begin{proposition}
\begin{equation}
F_L(2\La_0,\La_k)=
\sum
q^{\vec{m}^t C^{-1}\vec{m}}
\prod_{i=1}^{n-1}\Mult{\ell_i + m_i}{m_i}{q},
\end{equation}
with $\vec{\ell}$ given by (\ref{ldef}) and (\ref{rdef})
and with the sum
taken over all $\vec{m}\in (\Znn)^{(n-1)}$ satisfying
$k+\sum_{i=1}^{n-1} i m_i \equiv 0$.
\end{proposition}
\section{The general character}
In this section we calculate $F_L(\La_0+\La_j,\La_k)$ for arbitrary
$j$ ($0\le j\le n-1$). As we have already mentioned in section 2,
we count the weights in the principal picture, so that any fermionic
form can be reduced to one of the above form.

\begin{figure}[hbt]
\centerline{\epsffile{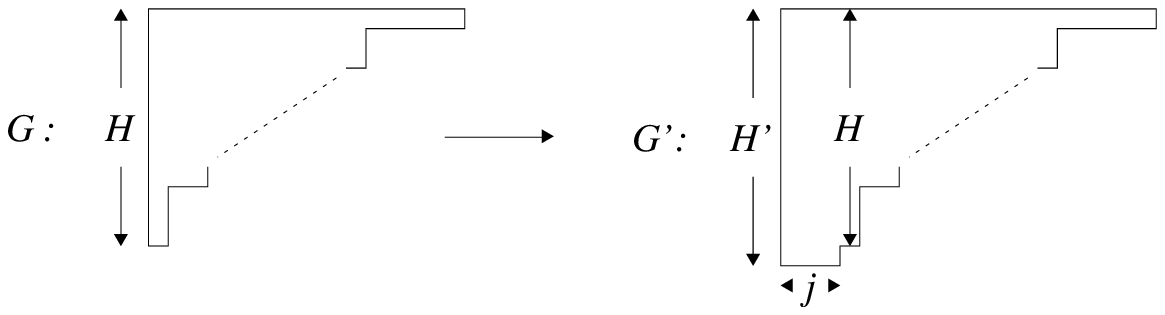}}
\caption{Embedding of a K-graph in
${\cal G}_L(\La_0+\La_j,\La_k)$  into ${\cal G}_{L+j}
(2\La_0,\La_k)$.}
\label{fig1sec8}
\end{figure}

First consider the following injection:
\begin{equation}
\begin{array}{rcl}
\P(\La_0+\La_j,\La_k) & \longrightarrow & {\cal P}_{L+j}(2\La_0,\La_k),\\
p=(\la_0,\cdots,\la_L) & \longmapsto & p',
\end{array}
\end{equation}
where $p'=(2\La_0,\La_0+\La_1,\cdots,\La_0+\La_{j-1},\la_0,\cdots,\la_L)$.
In terms of K-graphs, we have
\begin{equation}
\begin{array}{rcl}
\G(\La_0+\La_j,\La_k) & \longrightarrow & {\cal G}_{L+j}(2\La_0,\La_k),\\
G & \longmapsto & G'.
\end{array}
\label{injection}
\end{equation}
$G'$ is obtained from $G$ by placing the rectangle of width $j$
and height $H'$ in the left hand side of $G$ (see Figure~\ref{fig1sec8}),
where, $H$ being the height of $G$, $H'$ is determined by $H'-H=0$ or $n-j$,
$H'\in n\Z-k$.
It is clear that the image under the injection (\ref{injection}) is
the set of K-graphs in ${\cal G}_{L+j}(2\La_0,\La_k)$ having the
lowest plain of width at least $j$.

Now let us recall that we have established the following bijection
in the preceding sections:
\begin{equation}
\P(2\La_0,\La_k) \quad \longleftrightarrow \quad
\{(\vec{m},(F_1,\cdots,F_{n-1}))\}.
\label{bijection}
\end{equation}
Here $\vec{m}^t=(m_1,\cdots,m_{n-1})$ satisfying (\ref{column1}),
characterizes the parent graph, and $F_i$ is the Ferrers graph of
a partition with largest part $\leq \ell_i$ and number of parts
$\leq m_i$.  Regarding the image
of the injection (\ref{injection}), the following question arises:
How can we characterize
K-graphs having the lowest plain of width at least $j$ as elements in the
right hand side of the bijection (\ref{bijection})?
The answer is given by
\begin{proposition}
A K-graph in
$\G(2\La_0,\La_k)$
has a lowest plain of width
$w_1 \geq j$,
iff,
for $n-j+1 \leq i \leq n-1$,
the smallest part in the Ferrers graph $F_i$,
has at least $i+j-n$ nodes.
\end{proposition}
To prove this,
let $G$ be such a K-graph, $P$ its parent graph, and $M(G)$ and $M(P)$
the corresponding interpolating matrices.
We write $M(G)$ as in (\ref{Mwh}) and $M(P)$ as
$$
M(P)=\left(
\begin{array}{ccccc}
c_0+c_1 & \cdots & c_{i-1}+c_{i} & \cdots & c_{N-1}+c_N \\
n-c_1 & \cdots & n-c_{i} & \cdots & n-c_{N}
\end{array}
\right)
\qquad(c_0=0).
$$
\begin{quote}
{\em After removing g-components such that the first $i$ columns
of $M(G)$ equal $M(P)$, we must have that $w_{i+1}\geq j+c_i$.}
\end{quote}
We wish to prove the above assertion by induction on $i$.
For $i=0$, the assertion is clear from the assumption of the
proposition.
Next, let us assume the assertion for $i-1$. Let us also
assume that we have arrived at the minimal gap corresponding
to $c_{i-1}+c_i$. In order to prove the assertion for
$i$, we have to show that we can remove
an $(n-c_i)$-component, $j-c_{i-1}$ times strictly
horizontally. This can indeed be shown through straightforward,
though tedious consideration of the profile, and the conditions
on its various segments.
Proposition 2 follows from the above statement.
\medskip

Applying proposition 2, we immediately get
\begin{equation}
F_L(\La_0+\La_j,\La_k)=\sum q^{\vec{m}^tC^{-1}\vec{m}+Q}
\prod_{i=1}^{n-1}\Mult{\ell'_i+m_i}{m_i}{q},
\end{equation}
where
\begin{eqnarray}
Q & = & \sum_{i=n-j+1}^{n-1}(i+j-n)m_i-\frac{j}{n}\sum_{i=1}^{n-1}im_i
=-\vec{m}^tC^{-1}\vec{\e}_{n-j}, \nonumber \\
\ell'_i & = & \ell_i|_{L\rightarrow L+j}-(i+j-n)\theta(i+j-n),\\
& & \nonumber \\
& = & \left[C^{-1}(L\vec{\e}_{n-1}
+\vec{\e}_r+\vec{\e}_{n-j}-2\vec{m})\right]_i. \nonumber
\end{eqnarray}
Here $\ell_i$ is defined in (\ref{ldef}),
and $\theta$ in (\ref{step}).

In conclusion, we have the general form of the fermionic sum:
\begin{theorem}\label{Tf}
\begin{eqnarray}
F_L(\La_0+\La_j,\La_k) & = &
\sum q^{\vec{m}^tC^{-1}\vec{m}-\vec{m}^tC^{-1}\vec{\e}_{n-j}}
\prod_{i=1}^{n-1}\Mult{\ell_i+m_i}{m_i}{q}, \nonumber \\
\vec{\ell} & = &
C^{-1}(L\vec{\e}_{n-1}+\vec{\e}_r+\vec{\e}_{n-j}-2\vec{m}),
\end{eqnarray}
where the
sum is taken over all $\vec{m}\in(\Znn)^{n-1}$ satisfying
$k+\sum_{i=1}^{n-1}i m_i\equiv 0$, and with $r$ determined from
$L+j-2k\equiv r$, $0< r\leq n$.
\end{theorem}

\section{Summary and discussion}
In this paper we have presented a method to
compute finite analogues of the branching functions
of the coset
\begin{equation}
\widehat{sl(n)}_1 \otimes \widehat{sl(n)}_1 /
\widehat{sl(n)}_2.
\label{coset}
\end{equation}
Our approach, based on a direct counting of
Ferrers graphs related to the crystal base
formulation of the HWM's of $\widehat{sl(n)}$, leads
to what are known as {\em fermionic polynomials}.
This complements earlier results of Ref.~\cite{JMO} where
the same finite analogues of branching functions were
computed, and the result was expressed in terms of
{\em bosonic polynomials}.

Equating these two results, as formulated in the
Theorems~\ref{Tb} and \ref{Tf}, using equation (\ref{BtoF}),
we obtain the main
result of this paper:
\begin{theorem}
Let $\vec{m}$ and $\vec{\e}_i$ be $(n-1)$-dimensional vectors
with entries $(\vec{m})_i=m_i$ and $(\vec{\e}_i)_j=\delta_{i,j}$,
respectively. Also, let $C$ be the Cartan matrix of $sl(n)$
and ${\cal W}$ the Weyl group of $\widehat{sl(n)}$.
Defining the function
$b_{L,i}$ as in (\ref{bdef}),
the following polynomial identity holds for all $j,k=0,\ldots,n-1$:
\begin{eqnarray}
\lefteqn{
\sum_{\vec{m}\in (\Zs_{\geq 0})^{n-1}}
q^{\vec{m}^t C^{-1}\vec{m}-\vec{m}^tC^{-1}\vec{\e}_{n-j}}
\prod_{i=1}^{n-1}\Mults{
\left(\vec{m}+C^{-1}(L\vec{\e}_{n-1}+\vec{\e}_r+\vec{\e}_{n-j}-2\vec{m})
\right)_i}{m_i}{q}} \nonumber \\
&=&
q^{(|\La_k|^2- |\La_j|^2)/2} \label{polyids} \\
& & \quad \times
\sum_{w\in {\cal W}} (\mbox{det} \; w) \;
b_{L,j-k}(\La_k+\La_{j-k+L}+\rho - w(\La_0+\La_j+\rho) ),
\nonumber
\end{eqnarray}
with $0<r\leq n$ fixed by $L+j-2k\equiv r$
and with the sum over $\vec{m}$ restricted by
\begin{equation}
\frac{k}{n} + \left( C^{-1} \vec{m} \right)_{n-1} \in \Z.
\label{sumrest}
\end{equation}
\end{theorem}

Letting $L\to \infty$ we obtain the following $q$-series
identities for the branching functions of the coset
(\ref{coset}).

\begin{corollary}
Let $Q$ the root lattice and $\overline{{\cal W}}$ the Weyl group
of $sl(n)$. Then
\begin{eqnarray}
\lefteqn{
q^{(|\La_j|^2- |\La_k|^2)/2}
\sum_{\vec{m}\in (\Zs_{\geq 0})^{n-1}}
\frac{
q^{\vec{m}^t C^{-1}\vec{m}-\vec{m}^tC^{-1}\vec{\e}_{n-j}}}
{\prod_{i=1}^{n-1} (q)_{m_i}}
} \nonumber \\
&=& \frac{q^{\textstyle
\frac{|\La_j+\rho|^2}{2(n+2)}-
\frac{|\La_k+\rho|^2}{2(n+1)}}}
{(q)_{\infty}^{n-1}} \label{qids} \\
& & \times
\sum_{w\in \overline{{\cal W}}} (\mbox{det} \; w) \;
\Theta_{(n+2)(\overline{\La}_k+\overline{\rho})
-(n+1)w(\overline{\La}_j+\overline{\rho})
,(n+1)(n+2)}(q), \nonumber
\end{eqnarray}
with the sum over $\vec{m}$ again restricted by (\ref{sumrest}),
and with $\Theta_{\lambda,\ell}$ defined by
\begin{equation}
\Theta_{\lambda,\ell}(q) = \sum_{\alpha\in Q}
q^{\ell|\alpha-\case{\lambda}{\ell}|^2/2},
\end{equation}
for $\lambda\in \sum_{i=1}^{n-1} \C \overline{\La}_i$.
\end{corollary}
We note that the left-hand side of (\ref{qids}) coincides with the
character expressions of Lepowsky and Primc~\cite{LP}
for the $\Z_n$-parafermion conformal field theory.

\vspace{1cm}

The polynomial identities (\ref{polyids}) proven in this
work are, strictly-speaking, not new, since under
level-rank duality they map onto identities related to the
coset
\begin{equation}
\widehat{sl(2)}_{n-1} \otimes \widehat{sl(2)}_1 /
\widehat{sl(2)}_n.
\end{equation}
The latter were conjectured in \cite{StonyBrook,Melzer},
and proven in \cite{Berkovich,W,FW}.
However, the proof presented here is {\em intrinsically}
of $\widehat{sl(n)}$ type, and we expect it admits
generalization to the more general coset
\begin{equation}
\widehat{sl(n)}_{\ell} \otimes \widehat{sl(n)}_m /
\widehat{sl(n)}_{\ell+m}.
\end{equation}
Results related to general $\widehat{sl(2)}$-type cosets
were discussed in \cite{NY}. The fermionic character
form for certain sectors of the higher-rank parafermions were
proven in~\cite{Georgiev}.

For the case of $\widehat{sl(2)}$, the paths considered
in this paper admit yet another representation in terms of
Ferrers graphs. These graphs, obeying entirely different
conditions than our K-graphs, were
introduced in \cite{ABBBFV}.
They are also more general, in the sense that they lead to
character expressions for all $\widehat{sl(2)}$ cosets of type
${\cal C}_{2, \ell, 1}$, including rational values of $\ell$.
Results for this type of cosets have been discussed in
\cite{FQ,BMcC}.

\subsection*{Note added:}

After this work was completed, it was brought to our attention
that the main concepts introduced in this work are analogous to,
though intriguingly different from, concepts that are essential
to the theory of modular representations of the symmetric group
\cite{JK}. In particular, our K-graphs of the coset
${\cal C}_{n, 1, 1}$,
are known as $n$-regular Young diagrams, our {\it parent graphs}
are analogous to $n$-{\it cores}, our Gaussian polynomials generate
the analogues of $n$-{\it quotients}, our g-components are analogous
to {\it hook-ribbons}, and our counting procedure is very much
related to the evaluation of Kostka-Green-Foulkes polynomials.
However, there are differences, due to the fact that our K-graphs
obey additional conditions.

Now, to make things even more intriguing, we also learned that
the conditions obeyed by our K-graphs are almost identical to,
though stronger than, those obeyed by Young diagrams that
parametrize irreducible representations of $sl(n)$ which remain
irreducible under restriction to $sl(n-1)$. We hope to report on
these interesting relationships in future publications.

\section*{Acknowledgement}

We wish to thank Ch. Bessenrodt, A. Lascoux, B. Leclerc,
M. van Leeuwen, T. Nakajima, D. Stanton, J.-Y. Thibon,
H.-F. Yamada, and M. Yang for stimulating discussions
from which we learnt about the connections with the modular
representations of the symmetric group, and for their interest
in this work.
M.O. would like to thank the Department of Mathematics,
The University of Melbourne, for hospitality during a
stay (March - April, 1995), when most of this work was done.
This research was supported by the Australian Research
Council. One of the authors (M.O.) is partially supported by
Grant-in-Aid for Scientific Research on Priority Areas,
the Ministry of Education, Science and Culture, Japan.


\begin{thebibliography}{99}

\bibitem{Kac}
V.~G.~Kac,
{\em Infinite dimensional Lie algebras}
(Birkh\"auser, Boston, 1983).

\bibitem{KMPS}
S.~Kass, R.~V.~Moody, J.~Patera and R.~Slansky,
{\em Affine Lie algebras, weight multiplicities,
and branching rules}, Vol.~{\bf I},
(University of California Press, Berkeley, 1990).

\bibitem{GO}
P. Goddard and D. Olive, eds.,
{\em Kac-Moody and Virasoro algebras},
(World Scientific, Singapore, 1988).

\bibitem{DJMO}
E.~Date, M.~Jimbo, T.~Miwa and M.~Okado,
Phys.\ Rev.\ {\bf B35}, 2105 (1987).

\bibitem{JMO}
M.~Jimbo, T.~Miwa and M.~Okado,
Nucl.\ Phys.\ {\bf B300} [FS22], 74 (1988).

\bibitem{FZ}
A.~B.~Zamolodchikov and V.~A.~Fateev,
Sov.\ Phys.\ JETP {\bf 62}, 215 (1985).

\bibitem{BPZ}
A.~A.~Belavin, A.~M.~Polyakov and A.~B.~Zamolodchikov,
J.\ Stat.\ Phys.\ {\bf 34}, 763 (1984);
Nucl.\ Phys.\ B {\bf 241}, 333 (1984).

\bibitem{LP}
J.~Lepowsky and M.~Primc,
{\em Structures of the standard modules for the affine
Lie algebra $A_1^{(1)}$} Contemporary Mathematics,
{\bf 46} (AMS, Providence,1985).

\bibitem{StonyBrook}
R.~Kedem, T.~R.~Klassen, B.~M.~McCoy and E.~Melzer,
Phys.\ Lett.\ B {\bf 304}, 263 (1993);
Phys.\ Lett.\ B {\bf 307}, 68 (1993).

\bibitem{BaxterBook}
R.~J.~Baxter,
{\em Exactly solved models in statistical mechanics}
(Academic Press, London, 1982).

\bibitem{Schur}
I.~Schur,
S.-B.\ Preuss.\ Akad.\ Wiss.\ Phys.-Math.\ Kl., 302 (1917).

\bibitem{AndrewsBook}
G.~E.~Andrews,
{\em The Theory of Partition}
(Addison-Wesley, Reading, Massachusetts, 1976).

\bibitem{JMMO}
M.~Jimbo, K.~Misra, T.~Miwa, and M.~Okado,
Comm.\ Math.\ Phys.\ {\bf 136}, 543 (1991).

\bibitem{DJO}
E.~Date, M.~Jimbo and M.~Okado,
Comm.\ Math.\ Phys.\ {\bf 155}, 47 (1993).

\bibitem{MM}
K.~C.~Misra and T.~Miwa,
Comm.\ Math.\ Phys.\ {\bf 134}, 79 (1990).

\bibitem{Melzer}
E.~Melzer,
Int.\ J.\ Mod.\ Phys.\ A {\bf 9}, 1115 (1994).

\bibitem{Berkovich}
A.~Berkovich,
Nucl.\ Phys.\ B {\bf 431}, 315 (1994).

\bibitem{W}
S.~O.~Warnaar,
{\em Fermionic solution of the Andrews--Baxter--Forrester model I:
unification of TBA and CTM methods},
preprint University of Melbourne No. 02-95, hep-th/9501134.
To appear in J.~Stat.~Phys. ; \newline
{\em Fermionic solution of the Andrews--Baxter--Forrester model II:
proof of Melzer's polynomial identities},
preprint University of Melbourne.

\bibitem{FW}
O.~Foda and S.~O.~Warnaar,
{\em A bijection which implies Melzer's polynomial identities:
the $\chi_{1,1}^{(p,p+1)}$ case},
preprint University of Melbourne No. 03-95, hep-th/9501088.
To appear in Lett.~Math.~Phys.

\bibitem{Georgiev}
G.~Georgiev,
{\em  Combinatorial constructions of modules for
infinite-dimensio\-nal  Lie algebras, I.  Principal
subspace}, preprint Rutgers University, hep-th/9412054;
\newline
{\em Combinatorial  constructions of modules for
infinite-dimensional Lie algebras, II. Parafermionic
space}, preprint Rutgers University, q-alg/9504024.

\bibitem{NY}
A.~Nakayashiki and Y.~Yamada,
{\em  Crystalline Spinon Basis for RSOS Models},
preprint Kyushu University, hep-th/9505083.

\bibitem{ABBBFV}
G.~E.~Andrews, R.~J.~Baxter, D.~M.~Bressoud, W.~J.~Burge,
P.~J.~Forrester and G.~Viennot,
Europ.\ J.\ Comb. {\bf 8}, 341 (1987).

\bibitem{FQ}
O.~Foda and Y.-H.~Quano,
{\em Virasoro character identities from the Andrews--Bailey
construction},
preprint University of Melbourne No. 26-94, hep-th/9408086.

\bibitem{BMcC}
A.~Berkovich and B.~M.~McCoy,
{\em Continued Fractions and Fermionic Representations for
Characters of $M(p,p')$ Minimal Models},
preprint BONN-TH-94-28, ITPSB 94-060, hep-th/9412030.

\bibitem{JK}
G.~James and A.~Kerber,
{\em The Representation Theory of the Symmetric Group}
(Addison-Wesley, Reading, Massachusetts, 1981).

\end{thebibliography}
\end{document}